\documentclass[12pt]{amsart}
\usepackage{fullpage,amssymb,amscd}

\theoremstyle{plain}

\newtheorem{theorem}{Theorem}

\newtheorem{proposition}{Proposition}[section]
\newtheorem{lemma}[proposition]{Lemma}

\theoremstyle{definition}

\theoremstyle{remark}

\newlength{\globalparindent}
\def\C{\mathbb C}
\def\d{\delta}
\def\g{\mathfrak g}
\def\L{\mathfrak l}
\def\l{\lambda}
\def\O{\mathcal O}
\def\t{\mathfrak t}
\def\R{\mathbb R}
\def\T{\mathcal T}
\def\tl{\tilde{\lambda}}
\def\trho{\tilde{\rho}}
\def\W{\mathfrak W}
\def\Z{\mathbb Z}
\def\DH{\operatorname{DH}}

\newcommand{\ad}{\operatorname{ad}}
\newcommand{\Ad}{\operatorname{Ad}}

\newcommand{\card}{\operatorname{Card}}
\newcommand{\diag}{\operatorname{diag}}

\newcommand{\N}{\operatorname{N}}
\newcommand{\Res}{\operatorname{Res}}
\newcommand{\rank}{\operatorname{rank}}
\newcommand{\tr}{\operatorname{tr}}
\newcommand{\volume}{\operatorname{Volume}}
\newcommand{\udot}{{\mathaccent\cdot\cup}}

\begin{document}

\title[Duistermaat-Heckman Measure for Coadjoint Orbits]
  {The Duistermaat-Heckman Measure for the Coadjoint Orbits of Compact 
   Semisimple Lie Groups}

\author{Ami~Haviv}
\address{Institute of Mathematics\\
        The Hebrew University\\
        Giv'at-Ram, Jerusalem 91904\\
        Israel}
\email{haviv@math.huji.ac.il}

\dedicatory{This is a draft. Your comments are welcome.}
\date{This edition: \today; \ \ First edition: December 1, 1997.}

\begin{abstract}
We apply the Guillemin-Lerman-Sternberg theorem to reprove a formula of
Heckman for the Duistermaat-Heckman measure associated to the coadjoint action
of $T$, a maximal torus of a compact semisimple Lie group $G$, on 
a regular coadjoint $G$-orbit in $\g^\star$, the dual space of the Lie algebra
of $G$. This formula is, in an appropriate sense, a limiting case of the 
Kostant multiplicity formula.
\end{abstract}

\maketitle

\tableofcontents

\section{Introduction}

Let $G$ be a compact connected semisimple Lie group and let $T$ be
a maximal torus of $G$.
In~\cite{Heckman:Orbits}, Heckman considered the asymptotic
behavior of the multiplicities of the representations of $T$ occuring
in the restriction to $T$ of an highest weight representation of $G$, as
given by the Kostant multiplicity formula\footnote{Actually, Heckman
works in a more general situation of restriction to any closed subgroup
$K$ of $G$. We treat only the case $K=T$.}. 
He obtained an `asymptotic multiplicity function' and, using an
integration formula of Harish-Chandra, proved this function is closely
related to the 
push forward of the Liouville measure of the coadjoint orbit passing 
through the highest weight of the representation. This push forward 
measure was later generalized  by Duistermaat and Heckman to any 
Hamiltonian torus action with a proper moment map.

Guillemin and Sternberg (\cite{GuilleminSternberg:GeomQuant}) were able
to derive Heckman's results in the framework of symplectic geometry. Later,
together with Lerman, they found a general formula for the 
Duistermaat-Heckman measure (\cite{GLS:SymplecticFibrations}). This 
formula of Guillemin-Lerman-Sternberg is
closely related to the exact stationary phase formula of 
Duistermaat-Heckman (\cite{DuistermaatHeckman:ESP}).

We apply the Guillemin-Lerman-Sternberg formula in the case of coadjoint 
orbits, and reprove that the asymptotic multiplicity function is
(the Radon-Nikodym derivative of) the Duistermaat-Heckman measure. The
idea of this computation was known to specialists, but it may have not  
been written up before.   
Some impetus to this approach is gained from the 
recent proof of the Guillemin-Lerman-Sternberg theorem via 
cobordism (\cite{GGK:Cobordism},~\cite{GGK:Moment},
~\cite{Karshon:NonCmpctCobord}).
 
{\bf Acknowledgment.}
This note was written as a final assignment for a course titled ``Group 
Actions on Manifolds'', that was given by Yael Karshon on spring 1997 at 
the Hebrew University of Jerusalem. I would like to thank her for suggesting
me the topic of this note and for several useful discussions. I also
thank E.~Lerman for his comments.  

\section{The Symplectic Structure of the Coadjoint Orbit} \label{sect:KKS}

Let $V$ be a real vector space. The tangent space of $V$ at 
a point $p$, $\T_p(V)$, is identified with $V$ via the map 
$\tau_p:V \rightarrow \T_p(V)$, which assigns to a vector $v \in V$ the
derivation in the direction of $v$.
If $M \subset V$ is an (immersed) submanifold and $p \in M$, we identify
$\T_p(M)$ with $W_p$, the subspace of $V$ which is the image of the composition
\[
 \T_p(M) \stackrel{i_\star}{\longrightarrow} \T_p(V) \stackrel{\tau_p^{-1}}{\longrightarrow} V \, , 
\]
where $i_\star$ is the differential of the inclusion map $i:M \rightarrow V$. 

Assume a Lie group $G$ acts {\em linearly} on $V$, such that the $G$-action
preserves $M$. The following proposition is easy to prove. 
\begin{proposition}
\label{prop:EquiIdent}
The identification of the tangent spaces of $M$ with subspaces of $V$ is
$G$-equivariant. That is, the following diagram is commutative:
\begin{equation*}
\begin{CD}
W_p           @>{\tau_p}>>             \T_p(M)      \\
@VV{g}V                            @VV{g_\star}V  \\
W_{g \ldotp p} @>{\tau_{g \ldotp p}}>>   \T_{g \ldotp p}(M)              
\end{CD}
\qquad p \in M, \, g \in G. 
\end{equation*}
\end{proposition}

Let $G$ be a compact connected Lie group with a Lie algebra $\g$. Recall that 
for an action of $G$ on a manifold $M$, the {\em generating vector field}
corresponding to an element $\xi \in \g$ is the vector field on $M$ whose 
value at a point $p \in M$ is
\begin{equation*}
 \xi_M(p)=\left.\frac{d}{dt}\right|_{t=0} (\exp(t\xi) \ldotp p) \in \T_p(M) \, . 
\end{equation*}

$G$ acts on $\g$ by the adjoint action, hence also on $\g^{\star}$,  the dual 
space of $\g$, by the coadjoint action:
\[
 g \longmapsto \Ad^{\star}(g^{-1})\, , \quad g \in G. 
\]
Fix $\l \in \g^{\star}$ and let $\O=G \ldotp \l$ be the orbit of $G$
through $\l$. $\O$ is an embedded submanifold of $\g^{\star}$ (by compactness 
of $G$). Its dimension is $\dim G-\dim G_\l$ ($G_\l$ is the stabilizer of 
$\l$).
\begin{proposition}
The generating vector field for the action of $G$ on $\O$, 
corresponding to an element $\xi \in \g$, is given by the formula 
\begin{equation}
\label{eqn:CoadjOrbTang}
\xi_\O(f)=-f([\xi, \cdot]) \, , \quad f \in \O.
\end{equation}
($-f([\xi, \cdot])$ is an element of $\g^{\star}$ upon replacing the dot 
with elements from $\g$.)
\end{proposition}
\begin{proof}
As discussed before, we identify $\T_f(\O)$ ($f \in \O$) with a linear
subspace of $\g^{\star}$. Then, for any $\eta \in \g$,
\begin{alignat*}{2}
 \left.\frac{d}{dt}\right|_{t=0} (\exp(t\xi) \ldotp f) (\eta)  
 & = \left.\frac{d}{dt}\right|_{t=0} f(\Ad(\exp(-t\xi))\eta) & & \\
 & =f(\left.\frac{d}{dt}\right|_{t=0}\Ad(\exp(-t\xi))\eta) & 
 &\quad \text{(by linearity of $f$)} \\
 & =-f([\xi, \eta]) \, . & &
\end{alignat*} 
\end{proof}

From general theory (see, for example,~\cite{Bryant:LieGroups}), 
we know that for each $f \in \O$ the map
\[
 \g \rightarrow \T_f(\O) \quad , \quad \xi \mapsto \xi_{\O}(f)=-f([\xi, \cdot]) 
\]
is onto, and its kernel is
\[
 \g_f=\{\xi \in \g \mid f([\xi, \cdot])=0\} \, , 
\]
the Lie algebra of $G_f$. Thus, $T_f(\O) \cong \g/\g_f$.  

If $\g$ is equipped with a $G$-invariant inner product\footnote{
Since $G$ is compact, such an invariant inner product always exists. In our
application, we can (and will) take it to be the Killing form of $\g$.}
, denoted $(\cdot,\cdot)$, then it follows that the orthogonal complement of 
$\g_f$ with respect to that inner product,
$\g_f^\perp$, is isomorphic to $\T_f(\O)$. This isomorphism is denoted
$\Psi_f$:
\[
 \Psi_f:\g_f^\perp \widetilde{\longrightarrow} \T_f(\O) \, , \quad f \in \O \, . 
\]
Moreover, by 
proposition~\ref{prop:EquiIdent} and the $G$-invariance of $(\cdot,\cdot)$, 
this  `subspace model' for the space tangent to $\O$ at $f$ is compatible 
with the $G$-actions on $\g$ and $\O$: the diagram
\begin{equation}
\label{eqn:EquiSubsp}
\begin{CD}
\g_f^\perp       @>{\Psi_f}>>  \T_f(\O)         \\
@VV{\Ad(g)}V                  @VV{\Ad^{\star}(g^{-1})}V  \\
\g_{g \ldotp f}^\perp  @>{\Psi_{g \ldotp f}}>> \T_{g \ldotp f}(\O) 
\end{CD}
\qquad\qquad f \in \O, \, g \in G,
\end{equation}
is commutative. We will use this model in the sequel.

The Kirillov-Kostant-Souriau $2$-form on $\O$ is defined by
\begin{equation*}
 \omega_f(\xi_\O(f), \eta_\O(f)) = -f([\xi, \eta]) \, , \quad f \in \O \, , \xi, \eta \in \g \, . 
\end{equation*}
The basic properties of $\omega$ are summarized in
\begin{proposition}
$\omega$ is a well-defined smooth non-degenerate closed $2$-form on $\O$ (so that
$(\O, \omega)$ is a symplectic manifold). Furthermore, $\omega$ is 
$G$-invariant, and the inclusion map $\Phi_G:\O \hookrightarrow \g^{\star}$ 
is a {\em moment map} (that is, $\Phi_G$ is $G$-equivariant and satisfies the 
equations
\[
 d\Phi_G^{\xi}=-\iota(\xi_O)\omega \, , \quad \xi \in \g \, , 
\]
where $\Phi_G^{\xi}=\langle\Phi_G, \xi \rangle:\O \rightarrow \R$ is 
the $\xi$-coordinate of $\Phi_G$).
\end{proposition}
The proof of this proposition may be found, for example, in 
~\cite[Lemma 7.22]{BGV:HeatKernel}. 

Note that since the manifold $\O$ has a symplectic structure, its dimension 
is even. We set $n=\frac{1}{2}\dim \O$.

From now on we assume that $G$ is {\em semisimple} --- its center is a 
discrete (hence finite) subgroup of $G$. Let $T=G_{\l}$ be the stabilizer
of $\l$. We assume that $\l$ is {\em regular}, which means that $T$ is a 
maximal torus in $G$. We denote the Lie algebra of $T$ by $\t$. 

The coadjoint action of $G$ restricts to an action
of $T$ on $(\O, \omega)$. This action has a moment map
\[
 \Phi_T:\O \rightarrow \t^{\star} \quad , \quad \Phi_T=i^{\star} \circ \Phi_G=i^{\star} \, , 
\]
where $i^{\star}$ is the (restriction to $\O$ of the) projection
$i^{\star}:\g^{\star} \rightarrow \t^{\star}$, which is dual to the 
inclusion map $i:\t \rightarrow \g$.

The {\em Duistermaat-Heckman measure} corresponding to the action of $T$
on the symplectic manifold $(\O,\omega)$ is the signed measure on $\t^{\star}$
defined by 
\[
 \DH_T(W)=\int_{\Phi_T^{-1}(W)} \omega^n/n!=\int_{(i^{\star})^{-1}(W)} \omega^n/n! \, , 
\]
for an open subset $W \subset \t^{\star}$ with compact closure.

\section{The Guillemin-Lerman-Sternberg Theorem} \label{sect:GLS}

The Guillemin-Lerman-Sternberg theorem 
(\cite[Theorem 3.3.3]{GLS:SymplecticFibrations})
provides a formula for $\DH_T$. We recall the statement of the
theorem, with slight changes of notation,  
from~\cite[Lecture 18]{Karshon:GAoM} (see also [GGK2]).

Let $T$ be a compact torus acting on a compact symplectic manifold 
$(M^{2n}, \omega)$, with a moment map 
$\Phi_T:M \rightarrow \t^{\star}$. We 
assume that the set of $T$-fixed points, $M^T$, is finite.

For each $p \in M^T$, $T$ acts on $\T_p(M)$ via the isotropy action.
There is a unique decomposition of $T_p(M)$ into direct sum of
irreducible $T$-representation. The trivial representation does not
occur in this decomposition (this follows from the finiteness of $M^T$
and the equivariant slice theorem), so each summand is $2$-dimensional. 
We recall how such irreducible $2$-dimensional representations of $T$ are
parameterized. Let
\[
 L = \ker(\exp: \t \rightarrow T) \, .  
\]
$L$ is a lattice in $\t$. Its dual lattice (the {\em weight lattice}) is
defined by
\[
 L^{\star} = \{\alpha \in \t^{\star} \mid \alpha(\xi) \in 2\pi{\mathbb Z} \,,\, \text{ for all }   \xi \in L \} \, . 
\]
For $\alpha \in L^{\star}-\{0\}$, $T$ acts irreducibly on $\R^2$ by
\begin{equation}
\label{eqn:T_RealReps}
 R_{\alpha}:\exp(\xi) \longmapsto \left(
   \begin{array}{rr} \cos\alpha(\xi) & -\sin\alpha(\xi) \\
                     \sin\alpha(\xi) &  \cos\alpha(\xi)
   \end{array} \right) \, , \quad \xi \in \t \, . 
\end{equation}
The representations corresponding to $\pm \alpha$ are equivalent (via 
conjugation by 
$\left( \begin{smallmatrix} 0 & 1 \\ 1 & 0 \end{smallmatrix} \right)\,$).
Thus, the irreducible $2$-dimensional representations of $T$ are 
parameterized by $(L^{\star} -\{0\})/\pm1$.

We can also attach to $\alpha \in L^{\star}$ a $1$-dimensional
{\em complex} representation of $T$, namely 
\begin{equation}
\label{eqn:T_CpxReps}
 \tilde{R}_{\alpha}:\exp(\xi) \longmapsto \exp(\sqrt{-1}\alpha(\xi)) \, , \quad \xi \in \t \, . 
\end{equation}
There is an $\R$-linear isomorphism $\R^2 \tilde{\rightarrow} \C$ (sending
$(1,0) \mapsto 1, \, (0,1) \mapsto \sqrt{-1}$),
which intertwines the representations $R_{\alpha}$ and $\tilde{R}_{\alpha}$. 
When composing this isomorphism  with 
the complex conjugation, we get an intertwiner between $R_{-\alpha}$ and 
$\tilde{R}_{\alpha}$. Moreover, an $\R$-linear isomorphism 
$\R^2 \tilde{\rightarrow} \C$ with these intertwining properties is unique
up to multiplication by a non-zero {\em real} scalar. It follows that 
choosing one element from the unordered pair $\pm\alpha$ is the same thing 
as endowing $\R^2$ with an invariant complex structure.
 
Returning to the isotropy action of $T$ on $\T_p(M)$, we denote by
$\pm \alpha_{p,j}\, , \, j=1,\ldots, n$ the parameters of the irreducible
summands. ($\pm \alpha_{p,j}$ are called the {\em isotropy weights} at $p$.)
We now pick a `polarizing vector' $\Lambda \in \t$ such that 
$\alpha_{p,j}(\Lambda) \neq 0$ for all $p \in M^T$ and for all $j$, 
and choose out
of each pair $\pm\alpha_{p,j}$ the one taking a positive value at $\Lambda$ 
(that is, we fix the notation so that $\alpha_{p,j}(\Lambda) > 0$). As 
explained before, the choices we have made determine $\R$-linear isomorphisms
\[
 \Theta_p:\T_p(M) \widetilde{\longrightarrow} \C^n \, , \quad p \in M^T \, . 
\]
$\Theta_p$ intertwines the isotropy action of $T$ on $\T_p(M)$ with the
$n$-dimensional complex representation of $T$
\[
 \exp(\xi) \longmapsto \diag(\exp(\sqrt{-1}\alpha_{p,1}(\xi)) , \ldots, \exp(\sqrt{-1}\alpha_{p,n}(\xi))) \, , \quad \xi \in \t \, . 
\] 
We endow $\T_p(M)$ with the orientation determined by the symplectic
structure of $M$, and give $\C^n$ the usual complex orientation. Set
\[
 \epsilon_p=\begin{cases} 
             +1 & \text{ if $\Theta_p$ preserves orientation,} \\
             -1 & \text{ otherwise.}
            \end{cases} 
\]

We can now state the Guillemin-Lerman-Sternberg theorem. 
\begin{theorem}
\label{thm:GLS}
The Duistermaat-Heckman measure $\DH_T$ is absolutely continuous with respect
to the Lesbegue measure on $\t^{\star}$, and its density function 
(Radon-Nikodym derivative, also called the Duistermaat-Heckman function) 
$\rho(x)$ is given by the formula
\[
 \rho(x)=\sum_{p \in M^T} \epsilon_p\rho_p(x) \, , \quad x \in \t^{\star} \, , 
\]
where 
\[
 \rho_p(x)=\volume \{(x_1, \ldots, x_n) \in \R_+^n \mid \Phi_T(p)+\sum_{j=1}^{n} x_{j}\alpha_{p,j}=x \} \, . 
\]
(The volume here is  $(n - \dim T)$-dimensional; we shall recall its precise 
definition in Section~\ref{sect:compareKMF} below.) 
\end{theorem}

\section{Generalities about Compact Semisimple Lie Groups}

Before going on, we need to collect several facts 
from the structure theory of compact semisimple Lie groups and their Lie 
algebras. (The proofs of these facts may be found, for example, 
in~\cite{FultonHarris:RepTheory}, Lecture 26 in particular.)

Let $G$ be a compact connected semisimple Lie group with a maximal torus $T$.
Let $\g_\C=\g \otimes_\R \C$ be the complexification of $\g$. It is a complex
semisimple Lie algebra. $\t_\C=\t \oplus \sqrt{-1}\t$ is a 
{\em Cartan subalgebra} of $\g_\C$. The adjoint action of $G$ on $\g$ 
(resp. of $T$ on $\t$) extends to an action on $\g_\C$ (resp. $\t_\C$). 
The coadjoint actions extend to actions on $\g_{\C}^{\star}$, 
$\t_{\C}^{\star}$, the complex duals of $\g_\C$, $\t_\C$.

The {\em Killing form}, defined by
\[
 (\xi, \eta)=\tr(\ad\xi \circ \ad\eta) \, , \quad \xi,\eta \in \g_\C \, , 
\]
is a non-degenerate bilinear form on $\g_\C$. It is {\em invariant}, that is,
\[
 ([\xi, \eta], \zeta)=(\xi,[\eta,\zeta]) \, , \quad \xi,\eta,\zeta \in \g_\C \, , 
\]
\[
 (g \ldotp \xi, g \ldotp \eta)=(\xi,\eta) \, , \quad g \in G, \, \xi,\eta \in \g_\C \, . 
\]
The restrictions of $(\cdot,\cdot)$ to $\t_\C$, $\g$, and $\t$ are 
non-degenerate as well; in fact, on $\g$ (hence also on $\t$) the restriction 
of $(\cdot,\cdot)$ is negative definite (so that its negation is an invariant 
inner product on $\g$). 

Define, for $\alpha \in \t_{\C}^{\star}$, 
\[
 (\g_{\C})_\alpha=\{\xi \in \g_\C \mid [\tau, \xi]=\alpha(\tau)\xi \text{ for all } \tau \in \t_\C \} \, . 
\]
The set of {\em roots of $(\g_\C,\t_\C)$} is 
$\Delta=\{\alpha \in \t_\C^{\star}-\{0\} \mid (\g_\C)_\alpha \neq 0 \}$. If 
$\alpha \in \Delta$ is a root, then $\dim_\C (\g_\C)_\alpha = 1$.
The values of $\alpha$ on $\t$ are pure imaginary, so that $\Delta \subset 
\sqrt{-1}\t^{\star}$; furthermore, $\Delta \subset \sqrt{-1}L^{\star}$.

$\Delta$ can be decomposed (in several ways) into a disjoint union  
\[
 \Delta=\Delta_+ \udot \Delta_- 
\]
such that $\Delta_+$ (the {\em positive roots}) and $\Delta_-$ 
(the {\em negative roots}) satisfy the following conditions:
\begin{enumerate}
\item $\Delta_-=\{-\alpha \mid \alpha \in \Delta_+ \}$.
\item If $\alpha,\beta \in \Delta_+$ and $\gamma=\alpha+\beta \in \Delta$,
then $\gamma \in \Delta_+$.
\end{enumerate} 
Any such a disjoint union decomposition can be specified by a vector
$\Lambda \in \sqrt{-1}\t$ such that 
\begin{equation}
\label{eqn:polarize+}
 \alpha \in \Delta_+ \Leftrightarrow \alpha(\Lambda) > 0 \, . 
\end{equation}
We fix a choice of a set $\Delta_+$ of positive roots. The cardinality 
of $\Delta_+$ is\footnote{Note that this 
agrees with the former meaning of the notation $n$, because
$\dim \O = \dim G - \dim T$.}  $n=\frac{1}{2}(\dim G -\dim T)$.
 
$\g_\C$ has a direct sum decomposition (the {\em Cartan decomposition})
\[
 \g_\C=\t_\C \oplus \bigoplus_{\alpha \in \Delta_+}((\g_\C)_\alpha \oplus (\g_\C)_{-\alpha}) \, . 
\]
For each $\alpha \in \Delta_+$, there is a triplet of vectors 
$X_\alpha, Y_\alpha, H_\alpha \in \g_\C$, such that:
\begin{enumerate}
\item $X_\alpha \in (\g_\C)_\alpha$, $Y_\alpha \in (\g_\C)_{-\alpha}$, and
$H_\alpha \in \sqrt{-1}\t$.
\item  The standard $\mathfrak{sl}_2\C$ commutation relations hold:
\begin{equation}
\label{eqn:sl2}
   [H_\alpha,X_\alpha]=2X_\alpha \, , \quad
   [H_\alpha,Y_\alpha]=-2Y_\alpha \, , \quad
   [X_\alpha,Y_\alpha]=H_\alpha \, . 
\end{equation}
\item $\alpha(H_\alpha)=2$ .
\end{enumerate}
The vectors $\{\sqrt{-1}H_\alpha \mid \alpha \in \Delta_+ \}$ span the Lie
algebra $\t$ (over $\R$). 

$\g$ inherits from $\g_\C$ the {\em real Cartan decomposition}
\[
 \g=\t \oplus \bigoplus_{\alpha \in \Delta_+} \L_\alpha \, ,
\]
where $\L_\alpha=\g \cap ((\g_\C)_\alpha \oplus (\g_\C)_{-\alpha})$. $\L_\alpha$
is a real plane, which is spanned (over $\R$) by the vectors
\[
 U_\alpha=X_\alpha-Y_\alpha \quad , \quad V_\alpha=\sqrt{-1}(X_\alpha+Y_\alpha) \, . 
\] 
 
Let $\N(T)=\{g \in G \mid g^{-1}Tg=T \}$ be the normalizer of $T$. The 
{\em Weyl group} $\W=\N(T)/T$ is finite. We fix a set $\hat{\W} \subset \N(T)$ 
of representatives for the elements of $\W$.
$\W$ acts naturally on $T$, $\t$, $\t_\C$, and their duals. The sign of
$w \in \W$ (denoted $(-1)^w$) is $\det(w:\t \rightarrow \t)$. The action
of $w \in \W$ on $\t$ permutes the vectors 
$\{\sqrt{-1}H_\alpha \mid \alpha \in \Delta \}$, and 
\begin{equation}
\label{eqn:signw}
(-1)^{w}= \text{ parity of } 
\card\{\alpha \in \Delta_+ \mid w\ldotp(\sqrt{-1}H_\alpha)=\sqrt{-1}H_\beta,
\text{ for some } \beta \in \Delta_- \} \, .  
\end{equation} 

\section{Computation for a Coadjoint Orbit} \label{sect:compute}

We work out the various ingredients of Theorem~\ref{thm:GLS} for $M=\O$, the
coadjoint orbit of $G$ through a regular element $\l \in \g^{\star}$, with
the action of $T=G_{\l}$.
 
A point $p=g \ldotp \l \in \O$ is fixed by $T$ if and only if
\[
 t \ldotp g \ldotp \l=g \ldotp \l \, \text{ for all } t \in T \, .  
\]
This equality says that $g^{-1}tg \in G_\l=T$, that is, $g \in \N(T)$. 
Similarly, $g_1, g_2 \in \N(T)$ yield the same orbit element $p$ precisely when
$g_1T=g_2T$. We have shown:
\begin{proposition}
The fixed points of $T$ in $\O$ are
\begin{equation*}
 \O^T=\{w \ldotp \l \mid w \in \hat{\W} \} \, . 
\end{equation*}
\end{proposition}  
Hence, the formula of Theorem~\ref{thm:GLS} in this case takes the form
\begin{equation}
\label{eqn:rhoOrbit}
 \rho_\O(x)=\sum_{w \in \hat{\W}} \epsilon_{w \ldotp \l}\rho_{w \ldotp \l}(x) \, , \quad x \in \t^{\star} \, , 
\end{equation}
We turn to the analysis of the isotropy actions of $T$ on 
$\T_{w \ldotp \l}(\O)$. Bearing in mind that $\g_{w \ldotp \l}=\t$, the 
commutative diagram (\ref{eqn:EquiSubsp}) shows that, for all $t \in T$,
\[
 \Ad t(\xi)=\Psi_{w \ldotp \l}^{-1} (t \ldotp \Psi_{w \ldotp \l}(\xi)) \, , \quad \xi \in \t^{\perp} \, . 
\]  
This means that all the isotropy actions are equivalent to the adjoint action 
of $T$ on the orthogonal complement of $\t$ in $\g$ (with respect to 
the Killing form of $\g$).
\begin{lemma}
The orthogonal complement of $\t$ in $\g$ is
\begin{equation}
\label{eqn:t_perp}
 \t^{\perp}=\bigoplus_{\alpha \in \Delta_+} \L_\alpha \, . 
\end{equation}
\end{lemma}
\begin{proof}
By dimension considerations, it is enough to prove that 
$\L_\alpha \subset \t^{\perp}$, for all $\alpha \in \Delta_+$. 
Pick $\xi_0 \in \t$ such that $\alpha(\xi_0) \neq 0$. For each 
$\xi \in \t$,
\begin{align*} 
\alpha(\xi_0)(\xi, X_\alpha) &= (\xi, [\xi_0, X_\alpha]) \\
                             &= ([\xi, \xi_0], X_\alpha) \qquad \text{(by invariance of $(\cdot,\cdot)\,$)} \\
                             &= 0 \, ,
\end{align*}
hence $(\xi, X_\alpha)=0$. Similarly, $\t$ is orthogonal to $Y_\alpha$ (inside
$\g_\C$), so that $\L_\alpha \subset \t$.
\end{proof}

Let us determine how $T$ acts on $\L_\alpha$. Put 
$\alpha'=-\sqrt{-1}\alpha \in \t^{\star}$. For $t \in T$, $t=\exp(\xi)$ for
some $\xi \in \t$, and then $\Ad t=\exp(\ad \xi)$. Since
\begin{align*} 
 \ad \xi \begin{pmatrix} U_\alpha \\ V_\alpha 
         \end{pmatrix}
 &= \alpha(\xi) \begin{pmatrix} (X_\alpha+Y_\alpha) \\
                                \sqrt{-1}(X_\alpha-Y_\alpha) 
                \end{pmatrix} \\
 &= \alpha'(\xi) \begin{pmatrix} V_\alpha \\ -U_\alpha 
                  \end{pmatrix} \, ,
\end{align*}
the action of $\Ad t$ on $\L_\alpha$ is represented (with respect to the basis
$(U_\alpha, V_\alpha) \, $) by the matrix 
\[
 \exp\left(\begin{smallmatrix} 0 & -\alpha'(\xi) \\ \alpha'(\xi) & 0
           \end{smallmatrix}\right) 
 =\left(\begin{smallmatrix}  \cos\alpha'(\xi) & -\sin\alpha'(\xi) \\
                             \sin\alpha'(\xi) &  \cos\alpha'(\xi)
        \end{smallmatrix}\right) \, .
\]
Thus, the $T$-action on $\L_\alpha$ is equivalent to the representation
$R_{\alpha'}$ of $T$ (defined by~\eqref{eqn:T_RealReps}). 
By \eqref{eqn:t_perp}, the isotropy weights (for each 
fixed point) are
$\{\pm \alpha' \mid \alpha \in \Delta_+ \}$. We can take 
$\Lambda'=\sqrt{-1}\Lambda \in \t$ as a polarizing vector, since by
(\ref{eqn:polarize+}),
\[
 \alpha'(\Lambda')=-\sqrt{-1}\alpha(\sqrt{-1}\Lambda)=\alpha(\Lambda) > 0 \, .
\]
We conclude:
\begin{proposition}
\label{prop:local_rho}
In~\eqref{eqn:rhoOrbit}, 
\[
\rho_{w \ldotp \l}(x)=\volume \{(x_\alpha)_{\alpha \in \Delta_+} \in \R_+^n \mid w \ldotp i^{\star}(\l) + \sum_{\alpha \in \Delta_+} x_{\alpha}\alpha'=x \} \, , \quad x \in \t^{\star} \, . 
\]
\end{proposition}

It remains to compute $\epsilon_{w \ldotp \l}$. It follows from the
 discussion preceding the proposition that the complex structure of 
$\L_\alpha$, determined by the polarization, is given by
\begin{equation}
\label{eqn:CpxStruct}
U_\alpha \longmapsto 1 \, , \quad V_\alpha \longmapsto \sqrt{-1} \, .
\end{equation}
Denote this map by $\Theta:\t^{\perp} \rightarrow \C^n$. We have to check whether 
the maps
\[
\Theta_{w \ldotp \l}: \T_{w \ldotp \l}(\O) \longrightarrow \C^{n} \, , \quad
\Theta_{w \ldotp \l}=\Theta \circ \Psi_{w \ldotp \l}^{-1} \qquad
(w \in \hat{\W}) 
\]
are orientation-preserving. By linear algebra, the definition of symplectic 
orientation, and~\eqref{eqn:CpxStruct}, this is the case if and only
if the {\em Pfaffian} of the $2n \times 2n$ matrix, which is made of 
the $2 \times 2$ blocks
\[
 \begin{pmatrix}
   \omega_{w\ldotp\l}(\Psi_{w\ldotp\l}(U_{\alpha_i}),
                      \Psi_{w\ldotp\l}(U_{\alpha_j}))
 & \omega_{w\ldotp\l}(\Psi_{w\ldotp\l}(U_{\alpha_i}),
                      \Psi_{w\ldotp\l}(V_{\alpha_j})) \\
   \omega_{w\ldotp\l}(\Psi_{w\ldotp\l}(V_{\alpha_i}),
                      \Psi_{w\ldotp\l}(U_{\alpha_j}))
 & \omega_{w\ldotp\l}(\Psi_{w\ldotp\l}(V_{\alpha_i}),
                      \Psi_{w\ldotp\l}(V_{\alpha_j}))
 \end{pmatrix}
 \, , \quad i,j=1, \ldots , n \, ,
\] 
is positive. (We fixed an arbitrary ordering $\alpha_1, \ldots , \alpha_n$
of the positive roots.)

We claim that only the diagonal blocks are non-zero. This is a consequence 
of the following simple remarks:
\begin{itemize}
\item
By the definition of $\omega$, a typical element in the matrix is
of the form 
\[ 
 -(w\ldotp\l)([U_{\alpha_i}, U_{\alpha_j}])
 =-\l(w^{-1}\ldotp[U_{\alpha_i}, U_{\alpha_j}]) \, .
\]
\item
For all 
$\alpha, \beta \in \Delta, \, [(\g_\C)_\alpha, (\g_\C)_\beta] \subset
(\g_\C)_{\alpha + \beta}$, so that the commutator may have non-zero
intersection with $\t_\C$ only if $\beta=-\alpha$. 
\item
Translating the previous remark to the real planes  $\L_\alpha$, we find that 
\[ 
 \text{ if } i \neq j \text{, then } 
   [\L_{\alpha_i}, \L_{\alpha_j}] \subset \t^{\perp} \, .  
\]
\item
Cartan decomposition and ~\eqref{eqn:t_perp} imply that 
$\t^{\perp} = [\t, \g]$. 
\item 
As $T=G_\l$, $\l$ is killed by the generating vector fields
$\{\xi_\O \mid \xi \in \t \}$. This means (using (\ref{eqn:CoadjOrbTang}))
that $\l([\t, \g])=0$.
\item
Conclusion: an element of a non-diagonal block belongs to
\[
 \l(w^{-1}\ldotp[\L_{\alpha_i}, \L_{\alpha_j}]) \subset
 \l(w^{-1}\ldotp\t^{\perp}) =
 \l(\t^{\perp}) = \l([\t, \g]) = \{0\} \, . 
\]
\end{itemize}

It remains to deal with the diagonal blocks. By the commutation 
relations (\ref{eqn:sl2}),
\[
  [U_{\alpha_i}, V_{\alpha_i}]=2\sqrt{-1}H_{\alpha_i} \, , \quad i=1, \ldots , n \, .
\]
 Hence, the diagonal blocks are of the form
\[
\begin{pmatrix} 0 & -2\l(w^{-1}\ldotp(\sqrt{-1}H_{\alpha})) \\
                2\l(w^{-1}\ldotp (\sqrt{-1}H_{\alpha})) & 0
\end{pmatrix} \, .
\]
The Pfaffian is the product of the upper-right block elements, that is
\[
   P_{w\ldotp\l}=\prod_{\alpha \in \Delta_+}(-2\l(w^{-1}\ldotp(\sqrt{-1}H_{\alpha}))) \, .
\]
Let us concentrate on the case when 
$-\sqrt{-1}(i^{\star}\l) \in \t_\C^{\star}$
is {\em strongly dominant}, that is the inequality
\[
 -\sqrt{-1}\l(H_\alpha) > 0 
\]
holds for all $\alpha \in \Delta_+$. (We also term $\l$ as strongly dominant 
in this case.) Using (\ref{eqn:signw}) we get:
\begin{proposition}
If $\l$ is strongly dominant, then $P_{w \ldotp \l}$ is non-zero and its sign
equals $(-1)^{w}$.
\end{proposition}

This finishes the proof of
\begin{theorem}
\label{thm:main-thm}
Assume $\l \in \g^{\star}$ is strongly dominant. The  
Duistermaat-Heckman function for $\O=G\ldotp\l$ is
\[
 \rho_{\O}(x)=\sum_{w \in \hat{\W}} (-1)^w \rho_{w\ldotp\l}(x) \, , \quad x \in \t^{\star} \, ,
\]
where $\rho_{w \ldotp \l}(x)$ is given by the formula in 
Proposition~\ref{prop:local_rho}.
\end{theorem}

\section{Comparison with the Kostant Multiplicity Formula}
\label{sect:compareKMF}

In the preceding sections we considered a fixed linear functional
$\l \in \g^{\star}$ and took $T=G_{\l}$.
We now fix a maximal torus $T$ of $G$ and consider linear functionals 
$\l \in \t^{\star}$. Any such $\l$
can be extended to an element $\tl \in \g^{\star}$ via the direct sum 
decomposition $\g = \t \oplus \t^{\perp}$. Clearly, $i^{\star}\tl=\l$. Also, 
$T \subset G_{\tl}$, with equality (meaning that $\tl$ is regular) if and only
if $\l(\sqrt{-1}H_{\alpha}) \neq 0$ for all $\alpha \in \Delta$. 

Assuming that $\tl$ is strongly dominant (recall that this means that
\[
 \tl(-\sqrt{-1}H_{\alpha}) = \l(-\sqrt{-1}H_{\alpha}) > 0 \, ,
\] 
for all  $\alpha \in \Delta_+$), we deduce from Theorem~\ref{thm:main-thm} 
(applied to $\tl$) the following formula:
\[
 \rho_{G \ldotp \tl}(\mu) = \sum_{w \in \hat{\W}} (-1)^w \volume \{(x_\alpha)_{\alpha \in \Delta_+} \in \R_+^n \mid \sum_{\alpha \in \Delta_+} x_{\alpha}\alpha'=\mu - w \ldotp \l \} \, , \quad \mu \in \t^{\star} \, . 
\]
Recall that $\alpha'=-\sqrt{-1}\alpha$ . Put
\[
 \l'=-\sqrt{-1}\l \, , \quad \mu'=-\sqrt{-1}\mu \, ,
\]
and
$\trho_{\l}(\mu)=\rho_{G \ldotp \tl}(\mu) \,$.
With this notation we have:
\begin{equation}
\label{eqn:contKMF}
 \trho_{\l}(\mu) = \sum_{w \in \W} (-1)^w \volume \{(x_\alpha)_{\alpha \in \Delta_+} \in \R_+^n \mid \sum_{\alpha \in \Delta_+} x_{\alpha}\alpha= w \ldotp \l' - \mu' \} \, , \quad \mu' \in \sqrt{-1}\t^{\star} \, .  
\end{equation}

We now specialize to the case when  
$\l'$ is a strongly dominant weight (in the `complex' sense), that is,
$\tl$ is strongly dominant, and $\l \in L^{\star}$. The 
fundamental fact of the representation theory of $G$ is that there
exists a unique irreducible complex representation, $V(\l')$, of $\g_\C$ whose 
highest weight is $\l'$, and the restriction of $V(\l')$ to $\g$
lifts to an irreducible  unitary representation (also denoted $V(\l')$) of $G$.   
The restriction of $V(\l')$ from $G$ to $T$ is the direct sum of
irreducible (1-dimensional, complex) representations $\tilde{R}_{\mu}$
of $T$ (defined by~\eqref{eqn:T_CpxReps}).  Each 
representation $\tilde{R}_{\mu}$ occurs with a (possibly zero)
multiplicity which is denoted by $m_{\l}(\mu)$:
\[
 \Res_T V(\l') = \bigoplus_{\mu \in L^{\star}} m_{\l}(\mu) \tilde{R}_{\mu} \, . 
\] 
We let $Y(\l')$ be the set of $\mu \in L^{\star}$ such that 
$\tilde{R}_{\mu}$ has non-zero multiplicity in $\Res_T V(\l')$. 

The Kostant multiplicity formula (see, for example,
~\cite[Theorem IX.6.3]{Simon:RepsCmpct}) asserts that, 
for $\mu \in Y(\l')$, 
\begin{equation}
\label{eqn:discreteKMF}
 m_{\l}(\mu) = \sum_{w \in \W} (-1)^w \card \{(x_\alpha)_{\alpha \in \Delta_+} \in \Z_+^n \mid \sum_{\alpha \in \Delta_+} x_{\alpha}\alpha= w \ldotp (\l' + \d) - (\mu' + \d) \} \, , 
\end{equation}
where $\l'=-\sqrt{-1}\l \,$, $\mu'=-\sqrt{-1}\mu \,$,  and $\d$ is half 
the sum of the positive roots.

There is an apparent similarity between the formulas~\eqref{eqn:contKMF}
and~\eqref{eqn:discreteKMF}. We follow Heckman in quantifying this similarity
and quote several definitions and results from
~\cite[Section 2]{Heckman:Orbits}.

Let $A$ be a finite set contained in an open half space of a finite 
dimensional real vector space $E$. We also assume that $A$ is contained
in a lattice of maximal rank in $E$. Put $m=\card (A),\, r=\rank (A)$. 
The {\em partition function} of $A$ is
\[
 p_A(x) = \card \{ (x_\alpha)_{\alpha \in A} \in \Z_+^m \mid \sum_{\alpha \in A} x_{\alpha}\alpha = x \} \, , \quad x \in E \, .
\]
The {\em asymptotic partition function} of $A$ is
\[
 P_A(x)= \volume \{ (x_\alpha)_{\alpha \in A} \in \R_+^m \mid \sum_{\alpha \in A} x_{\alpha}\alpha = x \} \, , \quad x \in E \, ,
\]  
where the volume function is defined as follows. Fix an ordering
$A=(\alpha_1, \ldots , \alpha_m)$ of $A$. Let 
$W_A:\R^m \rightarrow E$ be the linear map sending $e_i$, the $i$-th standard 
basis vector of $\R^m$, to $\alpha_i$. $W_A^{-1}(0)$, the kernel of $W_A$, is 
an $(m-r)$-dimensional subspace of $\R^m$. By an easy argument, the integral 
points in $W_A^{-1}(0)$ form a lattice of maximal rank. Normalize the
Lesbegue measure on $W_A^{-1}(0)$ so that the fundamental
cell of this lattice has measure one. Call this normalized measure 
$\nu_0$. For every $x \in E$, $W_A^{-1}(x)$ is a translation of 
$W_A^{-1}(0)$, so we can translate the measure $\nu_0$ to get a measure
$\nu_x$ on $W_A^{-1}(x)$. (The translation invariance of Lesbegue measure
guarantees that there is no ambiguity in the definition of $\nu_x$.) Since
$A$ is contained in an open half space, $W_A^{-1}(x) \cap \R_+^m$ is a compact
convex polytope. The $\nu_x$-measure of this polytope is $P_A(x)$.

\begin{lemma}[Lemma 2.4 in~{\cite{Heckman:Orbits}}]
Suppose $\rank (A-\{\alpha\})=\rank (A)$ for all $\alpha \in A$. Fix $y$ in the
$\Z$-span of $A$. Then there exists a constant $C > 0$, depending on $y$, such 
that for all $x$ in the $\Z$-span of $A$
\begin{equation}
\label{eqn:estimate}
 |p_A(x+y)-P_A(x)| \leq C(1+|x|)^{m-r-1} \, .
\end{equation}
\end{lemma}

We want to apply this lemma with $E=\sqrt{-1}\t^{\star}$, $A=\Delta_+$,
$m=n(=\frac{1}{2}(\dim G - \dim T))$, $r=\dim T$,
$x_w = w \ldotp \l' - \mu'$ (with $\mu' \in Y(\l')$), and 
$y_w = w \ldotp \d - \d$ (for a fixed $w \in \W$).
Note that the $\Z$-span of $\Delta_+$ is the {\em root lattice} of $G$.
Let us examine the validity of the assumptions in the lemma:
\begin{itemize}
\item The condition $\rank(\Delta_+ - \{\alpha\})=\rank(\Delta_+)$, in the
case that $G$ is simple, merely rules out the possibility
$\g=\mathfrak{su}_2$. A similar restriction applies when $G$ is semisimple
(no $\mathfrak{su}_2$ direct factors). Since the result we wish to prove 
(Theorem~\ref{thm:MultLimit} below) holds rather trivially in this case,
this poses no serious difficulty.
\item Although $\d$ may not be a weight, $y_w$ does belong to the root lattice.
In fact,  it is the sum of distinct negative roots. (This can be easily 
deduced from, for example, the proof of Theorem IX.2.7 
in~\cite{Simon:RepsCmpct}.)  
\item The set $Y(\l')$ is invariant under the action of the Weyl group and 
the differences between its elements and $\l'$ lie in the root lattice 
(\cite[Section IX.4]{Simon:RepsCmpct}). These properties imply that
$x_w$  belongs to the root lattice. 
\end{itemize}

Thus, there is a constant $C_w$ such that the estimate (\ref{eqn:estimate}) 
holds. Summing over $w \in \W$, and using the obvious identities
\begin{align*}
 m_{\l}(\mu) &= \sum_{w \in \W} (-1)^w p_A(x_w+y_w) \, ,  \\
 \trho_{\l}(\mu) &= \sum_{w \in \W} (-1)^w P_A(x_w) \, ,  
\end{align*}
we infer that
\[
 |m_{\l}(\mu) - \trho_{\l}(\mu)| \leq C \sum_{w \in \W} (1+|w \ldotp \l' - \mu'|)^{n- \dim T - 1}
\]
(with $C = \max_{w \in \W}(C_w)$).

Set $s=n- \dim T$. The function $\trho_{\l}(\mu)$ has a simple 
homogeneity property:
\[
 \trho_{(k\l)}(k\mu) = k^s \trho_{\l}(\mu) \, , \quad k \in \Z_+ \, . 
\]
Now, for a constant $D$ depending on $\l$, $\mu$ (but not on $k$), we have
\[
 |m_{(k\l)}(k\mu) - \trho_{(k\l)}(k\mu)| \leq C \card(\W) (1+kD)^{s-1} \, ,
\]
so that
\[
k^s \left|\frac{m_{(k\l)}(k\mu)}{k^s} - \trho_{\l}(\mu)\right| \leq D_1 k^{s-1} \, ,
\] 
for another constant $D_1$. Dividing by $k^s$ and sending $k$ to infinity, 
we obtain
\begin{theorem}
\label{thm:MultLimit}
If $\l'$ is a strongly dominant weight and $\mu \in Y(\l')$, then
\[
 \trho_{\l}(\mu) = \lim_{k \rightarrow \infty} \frac{m_{(k\l)}(k\mu)}{k^s} \, .
\]
\end{theorem}

\ifx\undefined\bysame
        \newcommand{\bysame}{\leavevmode\hbox to3em{\hrulefill}\,}
\fi

\end{document}